\documentclass[12pt,a4paper]{article}
\usepackage[margin=2.5cm]{geometry}
\usepackage{color} 
\usepackage{xcolor} 
\usepackage{relsize} % relative size control
\usepackage{dsfont} %for identity operator
\usepackage{bbold} %for zero operator
\usepackage{amssymb}
\usepackage{amsmath}
\usepackage{graphicx}
\usepackage[cmtip,arrow]{xy}
\usepackage{pb-diagram,pb-xy}
\usepackage[noend]{algorithmic}
\usepackage{cancel}
\usepackage{hyperref}

\def\cmath{\color{blue}}
\def\ctxt{\color{black}}
\ifx\alert\undefined\newcommand{\alert}[1]{{\color{red}#1}}\fi
\newcommand{\baseform}[1]{\mathcal{A}_{#1}} % Basis permutation form
\def\baseformN{\mathrm{R}} % Number of base forms = rank of G
\def\bmat{\begin{pmatrix}}
\def\emat{\end{pmatrix}}

\newcommand{\bornform}[1]{\mathcal{B}_{#1}} % Born  form
\def\C{\mathbb{C}} % Complex numbers
\def\Q{\mathbb{Q}} % Рациональные числа
\def\R{\mathbb{R}} % Вещественные числа
\newcommand{\cabs}[1]{\left|#1\right|} % |#1|
\newcommand{\CyclG}[1]{\mathsf{C}_{#1}} % Cyclic group

\DeclareMathOperator{\diag}{diag}
\ifx\e\undefined\DeclareMathOperator{\e}{e}\fi

\newcommand{\farg}[1]{\!\left(#1\right)} % function arguments
\def\Hspace{\mathcal{H}} % Hilbert space
\def\imu{\mathrm{\mathbf{i}}} % imaginary unit
\DeclareMathOperator{\idmat}{\mathds{1}}
\DeclareMathOperator{\zeromat}{\mathbb{0}}
\newcommand{\IrrRep}[1]{\mathbf{#1}} % Неприводимое представление
\newcommand{\Math}[1]{$\cmath{}#1$} % Формулы с математическим цветом
\newcommand{\MathEq}[1]{\begin{equation*}\cmath{#1}\end{equation*}}
\newcommand{\MathEqLab}[2]{\begin{equation}\cmath{#1}\label{#2}\end{equation}}
\newcommand{\MathEqArr}[1]{\begin{align*}\cmath{}#1\end{align*}}

\def\NF{\mathcal{F}} % Любая числовая система

\newcommand{\ordset}[1]{\left[#1\right]} % Set of elements
\newcommand{\PermRep}[1]{\mathbf{\underline{#1}}} % Permutation representation
\def\regrep{\mathrm{P}} % Регулярное представление
\def\repq{\mathsf{U}} % Quantum representation symbol
\newcommand{\RNum}[1]{\uppercase\expandafter{\romannumeral #1\relax}} % Римские цифры
\newcommand{\set}[1]{\left\{#1\right\}} % Set of elements

\ifx\tr\undefined\DeclareMathOperator{\tr}{tr}\fi

\newcommand{\vect}[1]{\left(#1\right)} % Ordered set of elements, round brackets
\def\wG{\mathsf{G}} % Whole symmetry group: \mathcal{W}\mathrm{G}
\def\wS{\Omega} % Whole set of states:\mathrm{S}
\def\wSN{\mathsf{N}} % Size of whole set of states:{N_\wS}
\title{An Algorithm to Decompose Permutation Representations of Finite Groups:\\ {\Large{}A Polynomial Algebra Approach}}
\author{Vladimir V Kornyak\\[2pt]
\textit{Laboratory of Information Technologies}\\
\textit{Joint Institute for Nuclear Research}\\
\textit{Dubna, Moscow region, 141980, Russia}\\[2pt] \href{mailto:vkornyak@gmail.com}{vkornyak@gmail.com}}
\begin{document}
\maketitle
\begin{abstract}
\noindent
We describe an algorithm for splitting permutation representations of finite group over fields of characteristic zero into irreducible components.
The algorithm is based on the fact that the components of the invariant inner product in invariant subspaces are operators of projection into these subspaces.
An important element of the algorithm is the calculation of Gr\"{o}bner bases of polynomial ideals.
A preliminary implementation of the algorithm splits representations up to dimensions of tens of thousands.
Some examples of computations are given in Appendix \ref{appendix}.
\end{abstract}
\section{Introduction}\label{intro}
In general, the problem of splitting a module over an associative algebra into irreducible submodules is quite nontrivial.
An overview of the algorithmic aspects of this problem can be found in \cite{Holt}.
We consider here a particular, but important from our point of view \cite{Kornyak1}, case of the problem.
\par
Let \Math{\wG} be a permutation group on the set \Math{\wS\cong\set{1,\ldots,\wSN}}.
We will denote the action of \Math{g\in\wG} on \Math{i\in\wS} by \Math{i^g}.
To avoid inessential technical complications, we assume that \Math{\wG} acts transitively on \Math{\wS}.
\par
A representation of \Math{\wG} in an \Math{\wSN}-dimensional vector space over a field \Math{\NF} by the matrices \Math{\regrep\farg{g}}  with the entries 
\Math{\regrep\farg{g}_{ij}=\delta_{i^{\mathlarger{g}}j}}, 
where \Math{\delta_{ij}} is the Kronecker delta, will be called a \emph{permutation representation}.
Our approach uses the concept of inner product, so we have to exclude vector spaces over finite fields --- there are no reasonable inner products for such spaces.
Thus, we assume that the permutation representation space is a Hilbert space \Math{\Hspace_\wSN}.
From a constructive point of view it is sufficient to assume that the base field \Math{\NF} is a \emph{minimal splitting field} of the group \Math{\wG}.
Such field is a subfield of an \Math{m}-th cyclotomic field, where \Math{m} is a certain divisor of the \emph{exponent} of the group \Math{\wG}.
The field \Math{\NF}, being an abelian extension of the field \Math{\Q}, is a constructive dense subfield of \Math{\R} or \Math{\C}.  
\par
An orbit of \Math{\wG} on the Cartesian square \Math{\wS\times\wS} is called an \emph{orbital}  \cite{Cameron}.
The number of orbitals, called the \emph{rank} of \Math{\wG}, will be denoted by \Math{\baseformN}.
Among the orbitals of a transitive group there is one \emph{diagonal} orbital, \Math{\Delta_1=\set{\vect{i,i}\mid i\in\wS}},
which will always be fixed as the first element in the list of orbitals: \Math{\set{\Delta_1,\ldots,\Delta_\baseformN}}.
For transitive action of \Math{\wG} there is a natural one-to-one correspondence between the orbitals of \Math{\wG} and the orbits of a point stabilizer \Math{\wG_i}:
\MathEq{\Delta\longleftrightarrow\Sigma_i=\set{j\in\wS\mid \vect{i,j}\in\Delta}.}
The \Math{\wG_i}-orbits are called \emph{suborbits} and their cardinalities will be called the \emph{suborbit lengths}. 
Note that \Math{\cabs{\Delta}=\wSN\cabs{\Sigma_i}}.
\par
The invariance condition for a bilinear form \Math{A} in the Hilbert space \Math{\Hspace_\wSN} can be written as the system of equations 
\Math{A=\regrep\farg{g}A\regrep\farg{g^{-1}},~g\in\wG.} %{invcond}
It is easy to verify that in terms of the entries the equations of this system have the form \Math{\vect{A}_{\displaystyle{}ij}=\vect{A}_{\displaystyle{}i^gj^g}.}
Thus, the basis of all invariant bilinear forms is in one-to-one correspondence with the set of orbitals:
with each orbital \Math{\Delta_r\in\set{\Delta_1,\ldots,\Delta_\baseformN}} is associated a \emph{basis matrix} \Math{\baseform{r}} of the size \Math{\wSN\times\wSN} with the entries
\MathEq{
\vect{\baseform{r}}_{\displaystyle{}ij} =
\begin{cases}
1, &\text{if~}\vect{i,j}\in\Delta_r\,,\\
0, &\text{if~}\vect{i,j}\notin\Delta_r\,.
\end{cases}
}
In particular, for the diagonal orbital we have \Math{\baseform{1}=\idmat_{\wSN}}.
The matrices 
\MathEqLab{\baseform{1},\baseform{2},\ldots,\baseform{\baseformN}}{basis}
form a basis of the \emph{centralizer algebra} of the representation \Math{\regrep}.
The multiplication table for this basis has the form
\MathEqLab{\baseform{p}\baseform{q}=\sum_{r=1}^{\baseformN}C_{pq}^r\baseform{r},}{multtab}
where \Math{C_{pq}^r} are non-negative integers.
The commutativity of the centralizer algebra indicates that the permutation representation \Math{\regrep} is \emph{multiplicity-free}.
\section{Algorithm}\label{algo}
Let \Math{T} be a transformation (we can assume that \Math{T} is a unitary matrix) that splits the permutation representation \Math{\regrep} into \Math{M} irreducible components:
\MathEq{T^{-1}\regrep\farg{g}T=1\oplus\repq_{\!d_2}\farg{g}\oplus\cdots\oplus\repq_{\!d_m}\farg{g}\oplus\cdots\oplus\repq_{\!d_M}\farg{g},}
where \Math{\repq_{\!d_m}} is a \Math{d_m}-dimensional irreducible subrepresentation, \Math{\oplus} denotes the direct sum of matrices, i.e., \Math{A\oplus{}B=\diag\farg{A,B}}.
\par
The identity matrix \Math{\idmat_\wSN} is the \emph{standard inner product} in any orthonormal basis.
In the splitting basis we have the following decomposition of the standard inner product
\MathEq{\idmat_\wSN=\idmat_{d_1=1}\oplus\cdots\oplus\idmat_{d_m}\oplus\cdots\oplus\idmat_{d_M}.}
The \emph{inverse image} of this decomposition in the original permutation basis has the form
\MathEq{\idmat_\wSN=\bornform{1}+\cdots+\bornform{m}+\cdots+\bornform{M},}
where \Math{\bornform{m}} is defined by the relation
\MathEqLab{T^{-1}\bornform{m}T=\zeromat_{1+d_2+\cdots+d_{m-1}}\oplus\idmat_{d_m}\oplus\zeromat_{d_{m+1}+\cdots+d_M}.}{brel}
The main idea of the algorithm is based on the fact that \Math{\bornform{m}} is a \emph{projector}, i.e., \Math{\bornform{m}^2=\bornform{m}}.
Thus, all \Math{\bornform{m}}'s can be obtained as solutions of the equation 
\MathEqLab{X^2-X=0}{ideq}
for the generic invariant form 
\MathEq{X=x_1\baseform{1}+\cdots+x_\baseformN\baseform{\baseformN}.}
Using the multiplication table \eqref{multtab}, we can write the left-hand side of \eqref{ideq} as a set of polynomials
\MathEqLab{E\farg{x_1,\ldots,x_\baseformN}=\set{E_1\farg{x_1,\ldots,x_\baseformN},\ldots,E_\baseformN\farg{x_1,\ldots,x_\baseformN}}}{idpoly}
and equation \eqref{ideq} takes the form
\MathEqLab{E\farg{x_1,\ldots,x_\baseformN}=0.}{ideqpoly}
(This is an abbreviation for the set of equations \Math{\set{E_1=0,\ldots,E_\baseformN=0}}.)
Each polynomial in \eqref{idpoly} has the structure \Math{E_r\farg{x_1,\ldots,x_\baseformN}=Q_r\farg{x_1,\ldots,x_\baseformN}-x_r},
where  \Math{Q_r\farg{x_1,\ldots,x_\baseformN}} is a homogeneous polynomial of degree \Math{2} in the indeterminates \Math{x_1,\ldots,x_\baseformN}.
\par
In the basis \eqref{basis} the projector \Math{\bornform{m}} can be represented  as \MathEq{\bornform{m}=b_{m,1}\baseform{1}+b_{m,2}\baseform{2}+\cdots+b_{m,\baseformN}\baseform{\baseformN},}
where the vector 
\MathEq{B_m=\ordset{b_{m,1},\ldots,b_{m,\baseformN}}} 
is a solution of equation \eqref{ideqpoly}.
Since only \Math{\baseform{1}} has nonzero diagonal elements, we have \MathEq{\tr\bornform{m}=b_{m,1}\wSN.}
On the other hand, relation \eqref{brel} shows that \MathEq{\tr\bornform{m}=d_{m}.}
This allows us to fix the coefficient \Math{b_{m,1}}: \MathEq{b_{m,1}=d_m/\wSN.}
Thus the only relevant values of \Math{x_1} in \eqref{ideqpoly} are \Math{d/\wSN} for some \Math{d}'s from the interval \Math{\ordset{1..\wSN-1}}.
Any relevant natural number \Math{d} is either irreducible dimension or sum of such dimensions.
Using the orthogonality condition for the projectors, \MathEq{\bornform{m}\bornform{m'}=0 \text{~\ctxt{}if~} m\neq{}m',} we can exclude the consideration of dimensions that are sums of irreducible ones.
\par
The main part of our algorithm is a loop over the possible dimensions \Math{d}.
The loop starts with \Math{d=1} and ends when the sum of irreducible dimensions becomes equal to \Math{\wSN}.
The current \Math{d} is processed as follows:
\begin{enumerate}
	\item Substitute \Math{d} into the set of polynomials: \Math{E\farg{x_1,x_2,\ldots,x_\baseformN}\rightarrow{}E\farg{d/\wSN,x_2,\ldots,x_\baseformN}.}
	\item Compute the Gr\"{o}bner basis (\Math{\mathrm{Gb}}) of the polynomial system \Math{E\farg{d/\wSN,x_2,\ldots,x_\baseformN}.} 
	\item If \Math{\mathrm{Gb}=\ordset{1}} then the system of equations \Math{E\farg{d/\wSN,x_2,\ldots,x_\baseformN}=0} is inconsistent.
	Go to the next dimension: \Math{d\rightarrow{}d+1.} 
	\item Otherwise compute the Hilbert dimension (\Math{\mathrm{Hd}}) of the Gr\"{o}bner basis \Math{\mathrm{Gb}}.
	\item If \Math{\mathrm{Hd}=0} then all irreducible components of the dimension \Math{d} are multiplicity-free.
	If there are \Math{k} different \Math{d}-dimensional irreducible components then for the system of polynomial equations \Math{\mathrm{Gb}=0} we obtain
the following set of solutions%
\footnote{The computation of the solutions is always algorithmically realizable, since the problem involves only polynomial equations with abelian Galois groups.}
\MathEqLab{\set{B_{m+1}=\ordset{d/\wSN,b_{m+1,2},\ldots,b_{m+1,\baseformN}},\ldots,B_{m+k}=\ordset{d/\wSN,b_{m+k,2},\ldots,b_{m+k,\baseformN}}}.}{sols0}
Here \Math{m} is the number of irreducible components constructed before.
Apply the procedure \texttt{ProcessSingleSolution} (described below) to all solutions in \eqref{sols0}.
	\item 
If \Math{\mathrm{Hd}>0} then we encounter an irreducible component with a multiplicity \Math{\mathrm{k}>1.} 
The corresponding component of the centralizer algebra has the form \Math{A\otimes\idmat_d}, 
where \Math{A} is an arbitrary \Math{k\times{}k} matrix, and \Math{\otimes} denotes the Kronecker product.
The idempotency condition \Math{\vect{A\otimes\idmat_d}^2=A\otimes\idmat_d} implies \MathEqLab{A^2-A=0.}{idemA}
The complete family of solutions%
\footnote{It is well known that any solution of equation \eqref{idemA} can be represented as \Math{A=Q^{-1}\vect{\idmat_r\oplus\zeromat_{k-r}}Q,}
where \Math{Q} is an arbitrary invertible \Math{k\times{}k} matrix, and \Math{0\leq{r}\leq{k}}.} 
of \eqref{idemA} is a manifold of dimension 
\Math{\left\lfloor{\mathrm{k}^2}/{2}\right\rfloor=\mathrm{Hd}.}
Application of Gr\"{o}bner basis technique for polynomial systems with parameters involves the cumbersome procedure of partitioning the parameter space.
To avoid these difficulties, we construct \Math{k} mutually orthogonal particular solutions of the form \eqref{sols0} and apply to each of them the procedure \texttt{ProcessSingleSolution}.
\end{enumerate}
The input of procedure \texttt{ProcessSingleSolution} is the current irreducible projector \Math{\bornform{m}.}
The procedure performs the following.
\begin{enumerate}
	\item 
Calculate the orthogonality condition	\Math{\bornform{m}X=0.}	
This is a system of linear equations
\MathEq{O_m\farg{x_1,\ldots,x_\baseformN}=0}
for the indeterminates \Math{x_1,\ldots,x_\baseformN} with the numerical coefficients \Math{b_{m,1},\ldots,x_{m,\baseformN}.}
	\item
Add the orthogonality relations to the idempotency polynomial set	\eqref{idpoly}
\MathEq{E\farg{x_1,\ldots,x_\baseformN}\longrightarrow{}E\farg{x_1,\ldots,x_\baseformN}\cup{}O_m\farg{x_1,\ldots,x_\baseformN}.}
This is done in order to exclude from further consideration the subspace of the projector \Math{\bornform{m}.}
	\item Add \Math{\bornform{m}} to the list of irreducible projectors \Math{\set{\bornform{1},\ldots,\bornform{m-1}}\longrightarrow\set{\bornform{1},\ldots,\bornform{m-1},\bornform{m}}.}
\end{enumerate}
\section{Implementation}\label{imple}
Our approach involves some widely used methods of polynomial computer algebra.
Therefore it is reasonable, at least for the preliminary experience, to take advantage of computer algebra systems with the developed tools for working with polynomials.
\par
Our current implementation is a program written in \textbf{C}.
The input data for the program is a set of permutations \Math{S=\set{s_1,\ldots,s_K}} of degree \Math{\wSN} that generates the group \Math{\wG.}
\par
The program 
\begin{enumerate}
	\item computes the basis of the centralizer algebra \eqref{basis} and its multiplication table \eqref{multtab},
	\item constructs the system of quadratic polynomials \eqref{idpoly} --- the left hand sides of the idempotency condition \Math{X^2-X=0},
	\item constructs the bilinear system \Math{BX} corresponding to the orthogonality condition \Math{BX=0} --- 
the indeterminates \Math{B=\ordset{b_1,\ldots,b_\baseformN}} of this system	are replaced by specific numerical values in the procedure \texttt{ProcessSingleSolution},
\item generates the code for processing the above constructed polynomial data by the computer algebra system \textbf{Maple.}
\end{enumerate}
\section*{Conclusion}\label{conc}
The algorithm described here is based on the use of methods of polynomial algebra, which are considered algorithmically difficult.
However, our approach leads to a small number of low-degree polynomials.
As can be seen in Appendix \ref{appendix} even a straightforward implementation of the approach can cope with rather large tasks.
The data presented in the appendix shows that the most restrictive parameter for the \textbf{Maple} part of the implementation is the rank of representations, i.e., the number of polynomial indeterminates.
There are the obvious ways to improve performance: (1) to write in \textbf{C} polynomial algebra algorithms specialized for the  problem under consideration, (2) to use a more efficient system of computer algebra, for example, \textbf{Magma}.
\paragraph{Acknowledgement.}
I am grateful to V.P. Gerdt and Yu.A. Blinkov for fruitful discussions and valuable advices.

\appendix
\section{Examples of computations}\label{appendix}
\begin{itemize}
	\item Generators of representations are taken from the section {``Sporadic groups''} of the \textsc{Atlas} \cite{atlas}.
	\item The results presented below assume the following ordering for the centralizer algebra basis matrices 
\MathEq{\baseform{1}=\idmat_,\underbrace{\baseform{2},\ldots,\baseform{k},}_{\text{symmetric matrices}}
\underbrace{\baseform{k+1},\baseform{k+2}=\mathcal{A}^{\mathrm{T}}_{k+1},\ldots,\baseform{\baseformN-1},\baseform{\baseformN}=\mathcal{A}^{\mathrm{T}}_{\baseformN-1}}_{\text{asymmetric matrices}}.}
The matrices within the first sublist are ordered by the rule: \Math{A<B} if \Math{i_A<i_B}, where \Math{i_X= \min\vect{i\mid(X)_{i1}=1}}.
The same rule is applied to the first elements of the pairs of asymmetric matrices.
	\item Representations are denoted by their dimensions in bold (possibly with some signs added to distinguish different representations of the same dimension).
Permutation representations are underlined.
Multiple subrepresentations are underbraced in the decompositions.
	\item 
All timing data were obtained on a PC with 3.30GHz Intel Core i3 2120 CPU.
\end{itemize}

\subsection{Mathieu group \Math{M_{22}}}
\paragraph{Main properties:} 
Order:
\Math{\cabs{M_{22}}=443520=2^7\cdot3^2\cdot5\cdot7\cdot11.}\\
%\Math{\Exp\farg{{M_{22}}}=9240=2^3\cdot3\cdot5\cdot7\cdot11.}
Schur multiplier: \Math{\mathrm{M}\farg{M_{22}}=\CyclG{12}.} 
Outer automorphisms: \Math{\mathrm{Out}\farg{M_{22}}=\CyclG{2}.}

\subsubsection{\Math{770}-dimensional representation of \Math{M_{22}}}
Rank: \Math{9}. Suborbit lengths: \Math{1, 96, 144, 72, 144, 9, 16, 144, 144}.
\MathEq{\PermRep{770}\cong\IrrRep{1}\oplus\IrrRep{21}\oplus\underbrace{\vect{\IrrRep{55}\oplus\IrrRep{55}}}_{}\oplus\,\IrrRep{99}\oplus\IrrRep{154}\oplus\IrrRep{385}}
\vspace*{-20pt}
{\cmath	 
\begin{align*}
	\bornform{\IrrRep{1}} =&~ \frac{1}{770}\sum_{k=1}^{9}\baseform{k}
\\
	\bornform{\IrrRep{21}} =&~ \frac{3}{110}\vect{\baseform{1}+\frac{1}{12}\baseform{2}+\frac{1}{12}\baseform{3}+\frac{1}{12}\baseform{4}-\frac{3}{8}\baseform{5}+\baseform{6}-\frac{3}{8}\baseform{7}+\frac{1}{12}\baseform{8}+\frac{1}{12}\baseform{9}}
	\\
	\bornform{\IrrRep{55}} =&~ \frac{1}{14}\left\{\baseform{1}+\frac{1}{4}\vect{1-\imu\frac{3}{\sqrt{7}}}\baseform{2}-\frac{1}{4}\vect{1-\imu\frac{1}{\sqrt{7}}}\baseform{3}-\frac{1}{4}\vect{1+\imu\frac{1}{\sqrt{7}}}\baseform{4}\right.\\
	&\left.\hspace*{24pt}+\frac{1}{8}\vect{1+\imu\frac{1}{\sqrt{7}}}\baseform{5}+\vect{1+\imu\frac{2}{\sqrt{7}}}\baseform{6}+\frac{1}{8}\vect{1+\imu\frac{9}{\sqrt{7}}}\baseform{7}\right\}
	\\
	\bornform{\IrrRep{55}}{'} =&~ \frac{1}{14}\left\{\baseform{1}-\frac{1}{12}\vect{2-\imu\frac{9}{\sqrt{7}}}\baseform{2}+\frac{1}{12}\vect{2-\imu\frac{3}{\sqrt{7}}}\baseform{3}+\frac{1}{36}\vect{16+\imu\frac{9}{\sqrt{7}}}\baseform{4}\right.\\
	&\left.\hspace*{15.2pt}+\frac{1}{72}\vect{4-\imu\frac{9}{\sqrt{7}}}\baseform{5}-\frac{1}{9}\vect{1+\imu\frac{18}{\sqrt{7}}}\baseform{6}-\frac{1}{8}\vect{4+\imu\frac{9}{\sqrt{7}}}\baseform{7}-\frac{5}{36}\baseform{8}-\frac{5}{36}\baseform{9}\right\}
	\\
	\bornform{\IrrRep{99}} =&~ \frac{9}{70}\vect{\baseform{1}-\frac{1}{24}\baseform{2}-\frac{1}{9}\baseform{3}+\frac{1}{6}\baseform{4}-\frac{1}{24}\baseform{5}-\frac{1}{9}\baseform{6}+\frac{3}{8}\baseform{7}+\frac{1}{36}\baseform{8}+\frac{1}{36}\baseform{9}}
\end{align*}\begin{align*}	
	\bornform{\IrrRep{154}} =&~ \frac{1}{5}\vect{\baseform{1}+\frac{1}{12}\baseform{2}+\frac{1}{12}\baseform{3}-\frac{1}{18}\baseform{4}-\frac{1}{36}\baseform{5}-\frac{1}{9}\baseform{6}+\frac{1}{4}\baseform{7}-\frac{1}{18}\baseform{8}-\frac{1}{18}\baseform{9}}\hspace*{40pt}~
	\\
	\bornform{\IrrRep{385}} =&~ \frac{1}{2}\vect{\baseform{1}-\frac{1}{24}\baseform{2}-\frac{1}{18}\baseform{4}+\frac{1}{72}\baseform{5}-\frac{1}{9}\baseform{6}-\frac{1}{8}\baseform{7}+\frac{1}{36}\baseform{8}+\frac{1}{36}\baseform{9}}
\end{align*}
}
Time \textbf{C}: \Math{<1} sec. Time \textbf{Maple}: \Math{5} sec.

\subsubsection{\Math{990}-dimensional representation of \Math{3.M_{22}}}
Rank: \Math{13}. Suborbit lengths: \Math{1, 168, 336, 42, 7, 168, 168, 42, 42, 7, 7, 1, 1}.
\MathEqArr{\PermRep{990}\cong&\cmath\IrrRep{1}\oplus\IrrRep{21_\alpha}\oplus\IrrRep{21_\beta}\oplus\overline{\IrrRep{21_\beta}}\oplus\IrrRep{55}\oplus\IrrRep{99_\alpha}\oplus\overline{\IrrRep{99_\alpha}}\oplus\IrrRep{99_\beta}\\
&\cmath\oplus\IrrRep{105_+}\oplus\overline{\IrrRep{105_+}}\oplus\IrrRep{105_-}\oplus\overline{\IrrRep{105_-}}\oplus\IrrRep{154}}
\vspace*{-20pt}
{\cmath	 
\begin{align*}
	\bornform{\IrrRep{1}} =&~ \frac{1}{990}\sum_{k=1}^{13}\baseform{k}
\\
	\bornform{\IrrRep{21_\alpha}} =&~ \frac{7}{330}\left(\baseform{1}-\frac{5}{28}\baseform{2}+\frac{3}{14}\baseform{3}+\frac{3}{14}\baseform{4}-\frac{4}{7}\baseform{5}-\frac{5}{28}\baseform{6}-\frac{5}{28}\baseform{7}+\frac{3}{14}\baseform{8}+\frac{3}{14}\baseform{9}\right.\\
	&\left.\hspace*{24pt}-\frac{4}{7}\baseform{10}-\frac{4}{7}\baseform{11}+\baseform{12}+\baseform{13}\right)
\\
	\bornform{\IrrRep{21_\beta}} =&~ \frac{7}{330}\left\{\baseform{1}+\frac{1}{7}\baseform{2}+\frac{3}{7}\baseform{4}+\frac{5}{7}\baseform{5}-\frac{1}{14}\vect{1-\imu\sqrt{3}}\baseform{6}
	-\frac{1}{14}\vect{1+\imu\sqrt{3}}\baseform{7}\right.\\
	&\hspace*{24pt}-\frac{3}{14}\vect{1-\imu\sqrt{3}}\baseform{8}-\frac{3}{14}\vect{1+\imu\sqrt{3}}\baseform{9}-\frac{5}{14}\vect{1-\imu\sqrt{3}}\baseform{10}\\
	&\left.\hspace*{24pt}-\frac{5}{14}\vect{1+\imu\sqrt{3}}\baseform{11}-\frac{1}{2}\vect{1+\imu\sqrt{3}}\baseform{12}-\frac{1}{2}\vect{1-\imu\sqrt{3}}\baseform{13}\right\}
\\
	\bornform{\IrrRep{55}} =&~ \frac{1}{18}\left(\baseform{1}-\frac{1}{28}\baseform{2}-\frac{1}{14}\baseform{3}+\frac{3}{14}\baseform{4}+\frac{4}{7}\baseform{5}-\frac{1}{28}\baseform{6}-\frac{1}{28}\baseform{7}+\frac{3}{14}\baseform{8}+\frac{3}{14}\baseform{9}\right.\\
	&\left.\hspace*{24pt}+\frac{4}{7}\baseform{10}+\frac{4}{7}\baseform{11}+\baseform{12}+\baseform{13}\right)
\\
	\bornform{\IrrRep{99_\alpha}} =&~ \frac{1}{10}\left\{\baseform{1}+\frac{1}{12}\baseform{2}-\frac{1}{6}\baseform{4}-\frac{1}{24}\vect{1-\imu\sqrt{3}}\baseform{6}-\frac{1}{24}\vect{1+\imu\sqrt{3}}\baseform{7}\right.\\
	&\left.\hspace*{24pt}+\frac{1}{12}\vect{1-\imu\sqrt{3}}\baseform{8}+\frac{1}{12}\vect{1+\imu\sqrt{3}}\baseform{9}-\frac{1}{2}\vect{1+\imu\sqrt{3}}\baseform{12}\right.\\
	&\left.\hspace*{24pt}-\frac{1}{2}\vect{1-\imu\sqrt{3}}\baseform{13}\right\}
\\
	\bornform{\IrrRep{99_\beta}} =&~ \frac{1}{10}\left(\baseform{1}+\frac{1}{21}\baseform{2}-\frac{1}{14}\baseform{3}+\frac{1}{21}\baseform{4}-\frac{3}{7}\baseform{5}+\frac{1}{21}\baseform{6}+\frac{1}{21}\baseform{7}+\frac{1}{21}\baseform{8}+\frac{1}{21}\baseform{9}\right.\\
	&\left.\hspace*{24pt}-\frac{3}{7}\baseform{10}-\frac{3}{7}\baseform{11}+\baseform{12}+\baseform{13}\right)
	\\
\end{align*}\begin{align*}	
	\bornform{\IrrRep{105_\pm}} =&~ \frac{7}{66}\left\{\baseform{1}-\frac{1}{56}\vect{3\pm\frac{\sqrt{33}}{3}}\baseform{2}+\frac{1}{28}\vect{1\mp\frac{\sqrt{33}}{3}}\baseform{4}-\frac{1}{14}\vect{1\mp\sqrt{33}}\baseform{5}\right.\\
	&\left.\hspace*{24pt}+\frac{1}{112}\ordset{3\pm\frac{\sqrt{33}}{3}+\imu\vect{\sqrt{11}\pm3\sqrt{3}}}\baseform{6}+\frac{1}{112}\ordset{3\pm\frac{\sqrt{33}}{3}-\imu\vect{\sqrt{11}\pm3\sqrt{3}}}\baseform{7}\right.\\
	&\left.\hspace*{24pt}-\frac{1}{56}\ordset{1\mp\frac{\sqrt{33}}{3}-\imu\vect{\sqrt{11}\mp\sqrt{3}}}\baseform{8}-\frac{1}{56}\ordset{1\mp\frac{\sqrt{33}}{3}+\imu\vect{\sqrt{11}\mp\sqrt{3}}}\baseform{9}\right.\\
	&\left.\hspace*{24pt}+\frac{1}{28}\ordset{1\mp\sqrt{33}-\imu\vect{3\sqrt{11}\mp\sqrt{3}}}\baseform{10}+\frac{1}{28}\ordset{1\mp\sqrt{33}+\imu\vect{3\sqrt{11}\mp\sqrt{3}}}\baseform{11}\right.\\
	&\left.\hspace*{24pt}-\frac{1}{2}\vect{1\mp\imu\sqrt{3}}\baseform{12}-\frac{1}{2}\vect{1\pm\imu\sqrt{3}}\baseform{13}\right\}
\\
	\bornform{\IrrRep{154}} =&~ \frac{7}{45}\left(\baseform{1}+\frac{1}{28}\baseform{3}-\frac{1}{7}\baseform{4}+\frac{1}{7}\baseform{5}-\frac{1}{7}\baseform{8}-\frac{1}{7}\baseform{9}+\frac{1}{7}\baseform{10}+\frac{1}{7}\baseform{11}+\baseform{12}+\baseform{13}\right)
\end{align*}
}
Time \textbf{C}: \Math{<1} sec. Time \textbf{Maple}: \Math{37} sec.

\subsubsection{\Math{2016}-dimensional representation of \Math{3.M_{22}}}
Rank: \Math{16}. Suborbit lengths: \Math{1, 55, 165, 330, 165, 66, 66, 66, 330, 330, 165, 165, 55, 55, 1, 1}.
\MathEqArr{\PermRep{2016}\cong&\cmath\IrrRep{1}\oplus\IrrRep{21_\alpha}\oplus\IrrRep{21_\beta}\oplus\overline{\IrrRep{21_\beta}}\oplus\IrrRep{55}\oplus\IrrRep{105_+}\oplus\overline{\IrrRep{105_+}}\oplus\IrrRep{105_-}\oplus\overline{\IrrRep{105_-}}\\
&\cmath\oplus\IrrRep{154}\oplus\IrrRep{210_\alpha}\oplus\IrrRep{210_\beta}\oplus\overline{\IrrRep{210_\beta}}\oplus\IrrRep{231_\alpha}\oplus\IrrRep{231_\beta}\oplus\overline{\IrrRep{231_\beta}}}
\vspace*{-20pt}
{\cmath	 
\begin{align*}
	\bornform{\IrrRep{1}} =&~ \frac{1}{2016}\sum_{k=1}^{16}\baseform{k}
\\
	\bornform{\IrrRep{21_\alpha}} =&~ \frac{1}{96}\left(\baseform{1}-\frac{5}{11}\baseform{2}+\frac{3}{11}\baseform{3}-\frac{1}{11}\baseform{4}+\frac{3}{11}\baseform{5}-\frac{1}{11}\baseform{6}-\frac{1}{11}\baseform{7}-\frac{1}{11}\baseform{8}
	\right.\\&\left.\hspace*{24pt}
	-\frac{1}{11}\baseform{9}-\frac{1}{11}\baseform{10}+\frac{3}{11}\baseform{11}+\frac{3}{11}\baseform{12}-\frac{5}{11}\baseform{13}-\frac{5}{11}\baseform{14}+\baseform{15}+\baseform{16}\right)
	\\
	\bornform{\IrrRep{21_\beta}} =&~ \frac{1}{96}\left\{\baseform{1}+\frac{4}{11}\baseform{2}+\frac{3}{11}\baseform{3}-\frac{1}{11}\baseform{4}-\frac{4}{11}\baseform{6}
	+\frac{2}{11}\vect{1+\imu\sqrt{3}}\baseform{7}+\frac{2}{11}\vect{1-\imu\sqrt{3}}\baseform{8}
	\right.\\&\left.\hspace*{24pt}
+\frac{1}{22}\vect{1+\imu\sqrt{3}}\baseform{9}+\frac{1}{22}\vect{1-\imu\sqrt{3}}\baseform{10}
-\frac{3}{22}\vect{1+\imu\sqrt{3}}\baseform{11}-\frac{3}{22}\vect{1-\imu\sqrt{3}}\baseform{12}
	\right.\\&\left.\hspace*{24pt}
-\frac{2}{11}\vect{1+\imu\sqrt{3}}\baseform{13}-\frac{2}{11}\vect{1-\imu\sqrt{3}}\baseform{14}
-\frac{1}{2}\vect{1+\imu\sqrt{3}}\baseform{15}-\frac{1}{2}\vect{1-\imu\sqrt{3}}\baseform{16}\right\}
\\
	\bornform{\IrrRep{55}} =&~ \frac{55}{2016}\left(\baseform{1}+\frac{13}{55}\baseform{2}-\frac{1}{55}\baseform{3}-\frac{1}{55}\baseform{4}+\frac{13}{55}\baseform{5}-\frac{3}{11}\baseform{6}-\frac{3}{11}\baseform{7}-\frac{3}{11}\baseform{8}
	\right.\\&\left.\hspace*{24pt}
	-\frac{1}{55}\baseform{9}-\frac{1}{55}\baseform{10}-\frac{1}{55}\baseform{11}-\frac{1}{55}\baseform{12}+\frac{13}{55}\baseform{13}+\frac{13}{55}\baseform{14}+\baseform{15}+\baseform{16}\right)
\end{align*}\begin{align*}
	\bornform{\IrrRep{105_\pm}} =&~ \frac{5}{96}\left\{\baseform{1}-\frac{2}{55}\vect{1\pm\sqrt{33}}\baseform{2}+\frac{1}{55}\vect{4\pm\frac{\sqrt{33}}{3}}\baseform{3}
	\right.\\&\left.\hspace*{24pt}
	+\frac{1}{110}\vect{1\mp\frac{\sqrt{33}}{3}}\baseform{4}+\frac{1}{22}\vect{3\mp\frac{\sqrt{33}}{3}}\baseform{6}
	\right.\\
	&\left.\hspace*{24pt}
	-\frac{1}{44}\ordset{3\mp\frac{\sqrt{33}}{3}-\imu\vect{\sqrt{11}\mp3\sqrt{3}}}\baseform{7}-\frac{1}{44}\ordset{3\mp\frac{\sqrt{33}}{3}+\imu\vect{\sqrt{11}\mp3\sqrt{3}}}\baseform{8}
	\right.\\
	&\left.\hspace*{24pt}
	-\frac{1}{220}\ordset{1\mp\frac{\sqrt{33}}{3}-\imu\vect{\sqrt{11}\mp\sqrt{3}}}\baseform{9}-\frac{1}{220}\ordset{1\mp\frac{\sqrt{33}}{3}+\imu\vect{\sqrt{11}\mp\sqrt{3}}}\baseform{10}
	\right.\\&\left.\hspace*{24pt}
	-\frac{1}{110}\ordset{4\pm\frac{\sqrt{33}}{3}+\imu\vect{\sqrt{11}\pm4\sqrt{3}}}\baseform{11}-\frac{1}{110}\ordset{4\pm\frac{\sqrt{33}}{3}-\imu\vect{\sqrt{11}\pm4\sqrt{3}}}\baseform{12}
	\right.\\&\left.\hspace*{24pt}
		+\frac{1}{55}\ordset{1\pm\sqrt{33}+\imu\vect{3\sqrt{11}\pm\sqrt{3}}}\baseform{13}+\frac{1}{55}\ordset{1\pm\sqrt{33}-\imu\vect{3\sqrt{11}\pm\sqrt{3}}}\baseform{14}
	\right.\\&\left.\hspace*{24pt}
	-\frac{1}{2}\vect{1\pm\imu\sqrt{3}}\baseform{15}-\frac{1}{2}\vect{1\mp\imu\sqrt{3}}\baseform{16}\right\}
\\
	\bornform{\IrrRep{154}} =&~ \frac{11}{144}\left(\baseform{1}+\frac{7}{55}\baseform{2}+\frac{3}{55}\baseform{3}-\frac{3}{55}\baseform{4}-\frac{1}{11}\baseform{5}+\frac{1}{11}\baseform{6}+\frac{1}{11}\baseform{7}+\frac{1}{11}\baseform{8}
	\right.\\&\left.\hspace*{24pt}
	-\frac{3}{55}\baseform{9}-\frac{3}{55}\baseform{10}+\frac{3}{55}\baseform{11}+\frac{3}{55}\baseform{12}+\frac{7}{55}\baseform{13}+\frac{7}{55}\baseform{14}+\baseform{15}+\baseform{16}\right)
\\
	\bornform{\IrrRep{210_\alpha}} =&~ \frac{5}{48}\left(\baseform{1}-\frac{3}{55}\baseform{2}+\frac{1}{165}\baseform{3}+\frac{7}{165}\baseform{4}-\frac{7}{55}\baseform{5}-\frac{1}{11}\baseform{6}-\frac{1}{11}\baseform{7}-\frac{1}{11}\baseform{8}
	\right.\\&\left.\hspace*{24pt}
	+\frac{7}{165}\baseform{9}+\frac{7}{165}\baseform{10}+\frac{1}{165}\baseform{11}+\frac{1}{165}\baseform{12}-\frac{3}{55}\baseform{13}-\frac{3}{55}\baseform{14}+\baseform{15}+\baseform{16}\right)\\
	\bornform{\IrrRep{210_\beta}} =&~ \frac{5}{48}\left\{\baseform{1}-\frac{1}{15}\baseform{3}-\frac{1}{15}\baseform{4}
+\frac{1}{30}\vect{1+\imu\sqrt{3}}\baseform{9}+\frac{1}{30}\vect{1-\imu\sqrt{3}}\baseform{10}
	\right.\\&\left.\hspace*{24pt}
+\frac{1}{30}\vect{1+\imu\sqrt{3}}\baseform{11}+\frac{1}{30}\vect{1-\imu\sqrt{3}}\baseform{12}
	\right.\\&\left.\hspace*{24pt}
-\frac{1}{2}\vect{1+\imu\sqrt{3}}\baseform{15}-\frac{1}{2}\vect{1-\imu\sqrt{3}}\baseform{16}\right\}
	\\
	\bornform{\IrrRep{231_\alpha}} =&~ \frac{11}{96}\left(\baseform{1}-\frac{3}{55}\baseform{2}-\frac{1}{15}\baseform{3}+\frac{1}{165}\baseform{4}+\frac{1}{11}\baseform{5}+\frac{1}{11}\baseform{6}+\frac{1}{11}\baseform{7}+\frac{1}{11}\baseform{8}
	\right.\\&\left.\hspace*{24pt}
	+\frac{1}{165}\baseform{9}+\frac{1}{165}\baseform{10}-\frac{1}{15}\baseform{11}-\frac{1}{15}\baseform{12}-\frac{3}{55}\baseform{13}-\frac{3}{55}\baseform{14}+\baseform{15}+\baseform{16}\right)
\\	
	\bornform{\IrrRep{231_\beta}} =&~ \frac{11}{96}\left\{\baseform{1}-\frac{1}{33}\baseform{3}+\frac{2}{33}\baseform{4}-\frac{1}{11}\baseform{6}
+\frac{1}{22}\vect{1+\imu\sqrt{3}}\baseform{7}+\frac{1}{22}\vect{1-\imu\sqrt{3}}\baseform{8}
\right.\\&\left.\hspace*{24pt}
-\frac{1}{33}\vect{1+\imu\sqrt{3}}\baseform{9}-\frac{1}{33}\vect{1-\imu\sqrt{3}}\baseform{10}
	\right.\\&\left.\hspace*{24pt}
+\frac{1}{66}\vect{1+\imu\sqrt{3}}\baseform{11}+\frac{1}{66}\vect{1-\imu\sqrt{3}}\baseform{12}
	\right.\\&\left.\hspace*{24pt}
-\frac{1}{2}\vect{1+\imu\sqrt{3}}\baseform{15}-\frac{1}{2}\vect{1-\imu\sqrt{3}}\baseform{16}\right\}
\end{align*}
}
Time \textbf{C}: \Math{3} sec. Time \textbf{Maple}: \Math{1} h \Math{10} min \Math{48} sec.
\subsection{Higman--Sims group \Math{HS}}
\paragraph{Main properties:} 
\Math{\cabs{HS}=44352000=2^9\cdot3^2\cdot5^3\cdot7\cdot11.}~~
\Math{\mathrm{M}\farg{HS}=\CyclG{2}.}~~ 
\Math{\mathrm{Out}\farg{HS}=\CyclG{2}.}
\subsubsection{\Math{11200}-dimensional representation of \Math{2.HS}}
Rank: \Math{16}.\\ Suborbit lengths: \Math{1, 990, 660, 792, 792, 132, 110, 660, 132, 1, 1320, 1320, 1980, 1980, 165, 165}.
\MathEqArr{\PermRep{11200}\cong&\cmath\IrrRep{1}\oplus\IrrRep{22}\oplus\IrrRep{56}\oplus\IrrRep{77}\oplus\IrrRep{154}\oplus\IrrRep{175}\oplus\IrrRep{176}\oplus\overline{\IrrRep{176}}
\oplus\IrrRep{616}\oplus\overline{\IrrRep{616}}\\
&\cmath\oplus\IrrRep{770}\oplus\IrrRep{825}\oplus\IrrRep{1056}\oplus\IrrRep{1980}\oplus\overline{\IrrRep{1980}}\oplus\IrrRep{2520}}
\vspace*{-20pt}
{\cmath	 
\begin{align*}
	\bornform{\IrrRep{1}} =&~ \frac{1}{11200}\sum_{k=1}^{16}\baseform{k}
\\
\bornform{\IrrRep{22}} =&~ \frac{11}{5600}\left(\baseform{1}+\frac{13}{33}\baseform{2}-\frac{7}{33}\baseform{3}+\frac{1}{11}\baseform{4}+\frac{1}{11}\baseform{5}+\frac{13}{33}\baseform{6}+\frac{1}{11}\baseform{7}-\frac{7}{33}\baseform{8}
\right.\\&\left.\hspace*{32pt}
+\frac{13}{33}\baseform{9}+\baseform{10}-\frac{7}{33}\baseform{11}-\frac{7}{33}\baseform{12}+\frac{1}{11}\baseform{13}+\frac{1}{11}\baseform{14}-\frac{17}{33}\baseform{15}-\frac{17}{33}\baseform{16}
\right)
\\
\bornform{\IrrRep{56}} =&~ \frac{1}{200}\vect{\baseform{1}+\frac{1}{4}\baseform{3}+\frac{1}{4}\baseform{4}-\frac{1}{4}\baseform{5}+\frac{1}{4}\baseform{6}-\frac{1}{4}\baseform{8}-\frac{1}{4}\baseform{9}-\baseform{10}}
\\
\bornform{\IrrRep{77}} =&~ \frac{11}{1600}\left(\baseform{1}+\frac{1}{11}\baseform{2}+\frac{17}{132}\baseform{3}-\frac{23}{132}\baseform{4}-\frac{23}{132}\baseform{5}+\frac{37}{132}\baseform{6}-\frac{4}{11}\baseform{7}+\frac{17}{132}\baseform{8}
\right.\\&\left.\hspace*{32pt}
+\frac{37}{132}\baseform{9}+\baseform{10}-\frac{2}{33}\baseform{11}-\frac{2}{33}\baseform{12}+\frac{1}{66}\baseform{13}+\frac{1}{66}\baseform{14}+\frac{8}{33}\baseform{15}+\frac{8}{33}\baseform{16}
\right)
\\
\bornform{\IrrRep{154}} =&~ \frac{11}{800}\left(\baseform{1}+\frac{3}{55}\baseform{2}+\frac{7}{55}\baseform{3}+\frac{1}{11}\baseform{4}+\frac{1}{11}\baseform{5}-\frac{1}{11}\baseform{6}-\frac{19}{55}\baseform{7}+\frac{7}{55}\baseform{8}
\right.\\&\left.\hspace*{32pt}
-\frac{1}{11}\baseform{9}+\baseform{10}-\frac{1}{55}\baseform{11}-\frac{1}{55}\baseform{12}-\frac{3}{55}\baseform{13}-\frac{3}{55}\baseform{14}-\frac{7}{55}\baseform{15}-\frac{7}{55}\baseform{16}
\right)
\\
\bornform{\IrrRep{175}} =&~ \frac{1}{64}\left(\baseform{1}+\frac{7}{55}\baseform{2}-\frac{1}{15}\baseform{3}+\frac{1}{33}\baseform{4}+\frac{1}{33}\baseform{5}+\frac{1}{33}\baseform{6}+\frac{7}{55}\baseform{7}-\frac{1}{15}\baseform{8}
\right.\\&\left.\hspace*{32pt}
+\frac{1}{33}\baseform{9}+\baseform{10}+\frac{1}{33}\baseform{11}+\frac{1}{33}\baseform{12}-\frac{1}{15}\baseform{13}-\frac{1}{15}\baseform{14}+\frac{37}{165}\baseform{15}+\frac{37}{165}\baseform{16}
\right)
\\
\bornform{\IrrRep{176}} =&~ \frac{11}{700}\left(\baseform{1}+\frac{2}{33}\baseform{3}-\frac{1}{11}\baseform{4}+\frac{1}{11}\baseform{5}+\frac{7}{33}\baseform{6}-\frac{2}{33}\baseform{8}-\frac{7}{33}\baseform{9}-\baseform{10}
\right.\\&\left.\hspace*{32pt}
+\imu\frac{1}{33}\baseform{11}-\imu\frac{1}{33}\baseform{12}+\imu\frac{2}{33}\baseform{13}-\imu\frac{2}{33}\baseform{14}+\imu\frac{7}{33}\baseform{15}-\imu\frac{7}{33}\baseform{16}
\right)
\\
\bornform{\IrrRep{616}} =&~ \frac{11}{200}\left(\baseform{1}-\frac{7}{132}\baseform{3}+\frac{1}{44}\baseform{4}-\frac{1}{44}\baseform{5}+\frac{13}{132}\baseform{6}+\frac{7}{132}\baseform{8}-\frac{13}{132}\baseform{9}-\baseform{10}
\right.\\&\left.\hspace*{32pt}
-\imu\frac{1}{66}\baseform{11}+\imu\frac{1}{66}\baseform{12}-\imu\frac{1}{33}\baseform{13}+\imu\frac{1}{33}\baseform{14}+\imu\frac{4}{33}\baseform{15}-\imu\frac{4}{33}\baseform{16}
\right)
\\
\bornform{\IrrRep{770}} =&~ \frac{11}{160}\left(\baseform{1}-\frac{1}{165}\baseform{2}-\frac{1}{60}\baseform{3}-\frac{1}{44}\baseform{4}-\frac{1}{44}\baseform{5}+\frac{13}{132}\baseform{6}-\frac{4}{55}\baseform{7}-\frac{1}{60}\baseform{8}
\right.\\&\left.\hspace*{32pt}
+\frac{13}{132}\baseform{9}+\baseform{10}+\frac{7}{165}\baseform{11}+\frac{7}{165}\baseform{12}-\frac{1}{110}\baseform{13}-\frac{1}{110}\baseform{14}-\frac{16}{165}\baseform{15}-\frac{16}{165}\baseform{16}
\right)
\end{align*}\begin{align*}
\bornform{\IrrRep{825}} =&~ \frac{33}{448}\left(\baseform{1}+\frac{13}{495}\baseform{2}+\frac{7}{220}\baseform{3}-\frac{13}{396}\baseform{4}-\frac{13}{396}\baseform{5}-\frac{1}{12}\baseform{6}+\frac{12}{55}\baseform{7}+\frac{7}{220}\baseform{8}
\right.\\&\left.\hspace*{32pt}
-\frac{1}{12}\baseform{9}+\baseform{10}-\frac{1}{990}\baseform{13}-\frac{1}{990}\baseform{14}-\frac{8}{165}\baseform{15}-\frac{8}{165}\baseform{16}
\right)
\\
\bornform{\IrrRep{1056}} =&~ \frac{33}{350}\left(\baseform{1}-\frac{23}{495}\baseform{2}+\frac{3}{220}\baseform{3}+\frac{1}{36}\baseform{4}+\frac{1}{36}\baseform{5}+\frac{13}{132}\baseform{6}+\frac{6}{55}\baseform{7}+\frac{3}{220}\baseform{8}
\right.\\&\left.\hspace*{32pt}
+\frac{13}{132}\baseform{9}+\baseform{10}-\frac{1}{55}\baseform{11}-\frac{1}{55}\baseform{12}-\frac{2}{495}\baseform{13}-\frac{2}{495}\baseform{14}+\frac{4}{165}\baseform{15}+\frac{4}{165}\baseform{16}
\right)
\\
\bornform{\IrrRep{1980}} =&~ \frac{99}{560}\left(\baseform{1}+\frac{1}{132}\baseform{3}-\frac{1}{396}\baseform{4}+\frac{1}{396}\baseform{5}-\frac{7}{132}\baseform{6}-\frac{1}{132}\baseform{8}+\frac{7}{132}\baseform{9}-\baseform{10}
\right.\\&\left.\hspace*{32pt}
-\imu\frac{1}{33}\baseform{11}+\imu\frac{1}{33}\baseform{12}+\imu\frac{1}{99}\baseform{13}-\imu\frac{1}{99}\baseform{14}
\right)
\\
\bornform{\IrrRep{2520}} =&~ \frac{9}{40}\left(\baseform{1}-\frac{1}{165}\baseform{2}-\frac{1}{60}\baseform{3}+\frac{1}{396}\baseform{4}+\frac{1}{396}\baseform{5}-\frac{7}{132}\baseform{6}-\frac{4}{55}\baseform{7}-\frac{1}{60}\baseform{8}
\right.\\&\left.\hspace*{32pt}
-\frac{7}{132}\baseform{9}+\baseform{10}-\frac{1}{330}\baseform{11}-\frac{1}{330}\baseform{12}+\frac{1}{90}\baseform{13}+\frac{1}{90}\baseform{14}+\frac{4}{165}\baseform{15}+\frac{4}{165}\baseform{16}
\right)
\end{align*}
}
Time \textbf{C}: \Math{1} min \Math{12} sec. Time \textbf{Maple}: \Math{1} h \Math{39} min \Math{6} sec.

\subsection{McLaughlin group \Math{McL}}
\paragraph{Main properties:} 
\Math{\cabs{McL}=898128000=2^7\cdot3^6\cdot5^3\cdot7\cdot11.}~~
\Math{\mathrm{M}\farg{McL}=\CyclG{3}.}\\
\Math{\mathrm{Out}\farg{McL}=\CyclG{2}.}\\[5pt]
The \textsc{Atlas} \cite{atlas} contains four \Math{22275}-dimensional permutation representations of \Math{McL}, which are distinguished by the letters \Math{a}, \Math{b}, \Math{c} and \Math{d}.
\par
The computation shows that the representations \Math{\PermRep{22275a}}, \Math{\PermRep{22275c}} and \Math{\PermRep{22275d}} have the same structure of decompositions. 
In addition, the representations \Math{\PermRep{22275c}} and \Math{\PermRep{22275d}} have identical expressions for the irreducible projectors.
\subsubsection{Representation \Math{\PermRep{22275a}}}
Rank: \Math{13}.\\ Suborbit lengths: \Math{1, 4032, 112, 210, 210, 3360, 1260, 210, 1120, 3360, 3360, 2520, 2520}.
\MathEq{\PermRep{22275a}\cong\IrrRep{1}\oplus\IrrRep{22}\oplus\underbrace{\vect{\IrrRep{252}\oplus\IrrRep{252}}}_{}\oplus\underbrace{\vect{\IrrRep{1750}\oplus\IrrRep{1750}}}_{}\oplus\,\IrrRep{3520}\oplus\IrrRep{5103}\oplus\IrrRep{9625}}
\vspace*{-20pt}
{\cmath	 
\begin{align*}
	\bornform{\IrrRep{1}} =&~ \frac{1}{22275}\sum_{k=1}^{13}\baseform{k}
\\
\bornform{\IrrRep{22}} =&~ \frac{2}{2025}\left(\baseform{1}-\frac{1}{4}\baseform{2}+\frac{13}{28}\baseform{3}+\frac{9}{14}\baseform{4}+\frac{2}{7}\baseform{5}+\frac{2}{7}\baseform{6}-\frac{3}{7}\baseform{7}
\right.\\&\left.\hspace*{32pt}
+\frac{2}{7}\baseform{8}-\frac{1}{14}\baseform{9}+\frac{3}{28}\baseform{10}+\frac{3}{28}\baseform{11}-\frac{1}{14}\baseform{12}-\frac{1}{14}\baseform{13}
\right)
\end{align*}\begin{align*}
\bornform{\IrrRep{252}} =&~ \frac{28}{2475}\left\{\baseform{1}+\frac{1}{56}\baseform{2}+\frac{1}{1616}\vect{791-\frac{25\sqrt{33}}{7}}\baseform{3}+\frac{1}{2828}\vect{1307-\frac{25\sqrt{33}}{3}}\baseform{4}
\right.\\&\left.\hspace*{32pt}
-\frac{5}{2828}\vect{29+\sqrt{33}}\baseform{5}-\frac{5}{2828}\vect{29+\sqrt{33}}\baseform{6}+\frac{5}{5656}\vect{9-\frac{13\sqrt{33}}{3}}\baseform{7}
\right.\\&\left.\hspace*{32pt}
-\frac{5}{2828}\vect{29+\sqrt{33}}\baseform{8}+\frac{1}{11312}\vect{383+55\sqrt{33}}\baseform{9}
\right.\\&\left.\hspace*{32pt}
-\frac{1}{1616}\vect{19-\frac{5\sqrt{33}}{3}}\baseform{10}-\frac{1}{1616}\vect{19-\frac{5\sqrt{33}}{3}}\baseform{11}
\right\}
\\
\bornform{\IrrRep{252}}{'} =&~ \frac{28}{2475}\left\{\baseform{1}+\frac{1}{56}\baseform{2}-\frac{1}{11312}\vect{2911-{25\sqrt{33}}}\baseform{3}+\frac{1}{2828}\vect{\frac{2151}{5}+\frac{25\sqrt{33}}{3}}\baseform{4}
\right.\\&\left.\hspace*{32pt}
+\frac{1}{2828}\vect{\frac{523}{5}+5\sqrt{33}}\baseform{5}+\frac{1}{2828}\vect{\frac{523}{5}+5\sqrt{33}}\baseform{6}+\frac{1}{808}\vect{\frac{343}{5}+\frac{65\sqrt{33}}{21}}\baseform{7}
\right.\\&\left.\hspace*{32pt}
+\frac{1}{2828}\vect{\frac{523}{5}+5\sqrt{33}}\baseform{8}+\frac{11}{11312}\vect{\frac{689}{5}-5\sqrt{33}}\baseform{9}
\right.\\&\left.\hspace*{32pt}
-\frac{1}{1616}\vect{\frac{69}{7}+\frac{5\sqrt{33}}{3}}\baseform{10}-\frac{1}{1616}\vect{\frac{69}{7}+\frac{5\sqrt{33}}{3}}\baseform{11}-\frac{3}{35}\baseform{12}-\frac{3}{35}\baseform{13}
\right\}
\\
\bornform{\IrrRep{1750}} =&~ \frac{70}{891}\left\{\baseform{1}-\frac{1}{7}\baseform{2}+\frac{11}{14}\baseform{3}-\frac{3}{7}\baseform{4}+\frac{1}{140}\vect{73-9\imu\sqrt{231}}\baseform{5}-\frac{2}{35}\baseform{6}+\frac{3}{14}\baseform{7}
\right.\\&\left.\hspace*{32pt}
+\frac{1}{140}\vect{73+9\imu\sqrt{231}}\baseform{8}+\frac{13}{70}\baseform{9}+\frac{1}{70}\baseform{12}+\frac{1}{70}\baseform{13}
\right\}
\\
\bornform{\IrrRep{1750}}{'} =&~ \frac{70}{891}\left\{\baseform{1}+\frac{1}{7}\baseform{2}-\frac{13}{14}\baseform{3}+\frac{17}{35}\baseform{4}-\frac{1}{28}\vect{11-\imu\frac{9\sqrt{231}}{5}}\baseform{5}+\frac{2}{35}\baseform{6}-\frac{17}{70}\baseform{7}
\right.\\&\left.\hspace*{32pt}
-\frac{1}{28}\vect{11+\imu\frac{9\sqrt{231}}{5}}\baseform{8}-\frac{11}{70}\baseform{9}-\frac{1}{35}\baseform{10}-\frac{1}{35}\baseform{11}+\frac{1}{70}\baseform{12}+\frac{1}{70}\baseform{13}
\right\}
\\
\bornform{\IrrRep{3520}} =&~ \frac{64}{405}\left(\baseform{1}-\frac{1}{64}\baseform{2}-\frac{11}{112}\baseform{3}+\frac{3}{70}\baseform{4}-\frac{1}{70}\baseform{5}-\frac{1}{70}\baseform{6}+\frac{3}{140}\baseform{7}
\right.\\&\left.\hspace*{32pt}
-\frac{1}{70}\baseform{8}-\frac{17}{1120}\baseform{9}+\frac{3}{224}\baseform{10}+\frac{3}{224}\baseform{11}+\frac{1}{280}\baseform{12}+\frac{1}{280}\baseform{13}
\right)
\\
\bornform{\IrrRep{5103}} =&~ \frac{63}{275}\left(\baseform{1}-\frac{1}{144}\baseform{2}+\frac{1}{28}\baseform{3}-\frac{1}{35}\baseform{4}-\frac{13}{210}\baseform{5}+\frac{1}{80}\baseform{6}+\frac{1}{180}\baseform{7}
\right.\\&\left.\hspace*{32pt}
-\frac{13}{210}\baseform{8}+\frac{9}{560}\baseform{9}-\frac{1}{168}\baseform{10}-\frac{1}{168}\baseform{11}+\frac{1}{180}\baseform{12}+\frac{1}{180}\baseform{13}
\right)
\\
\bornform{\IrrRep{9625}} =&~ \frac{35}{81}\left(\baseform{1}+\frac{1}{112}\baseform{2}+\frac{1}{28}\baseform{3}-\frac{1}{35}\baseform{4}+\frac{1}{70}\baseform{5}-\frac{1}{560}\baseform{6}-\frac{1}{140}\baseform{7}
\right.\\&\left.\hspace*{32pt}
+\frac{1}{70}\baseform{8}-\frac{1}{80}\baseform{9}+\frac{1}{280}\baseform{10}+\frac{1}{280}\baseform{11}-\frac{1}{140}\baseform{12}-\frac{1}{140}\baseform{13}
\right)
\end{align*}
}
Time \textbf{C}: \Math{3} min \Math{32} sec. Time \textbf{Maple}: \Math{14} sec.
\subsubsection{Representation \Math{\PermRep{22275b}}}
Rank: \Math{6}. Suborbit lengths: \Math{1, 2240, 6720, 8064, 210, 5040}.
\MathEq{\PermRep{22275b}\cong\IrrRep{1}\oplus\IrrRep{252}\oplus\IrrRep{1750}\oplus\IrrRep{5103}\oplus\IrrRep{5544}\oplus\IrrRep{9625}}
\vspace*{-20pt}
{\cmath	 
\begin{align*}
	\bornform{\IrrRep{1}} =&~ \frac{1}{22275}\sum_{k=1}^{6}\baseform{k}
\\
\bornform{\IrrRep{252}} =&~ \frac{28}{2475}\left(\baseform{1}+\frac{13}{112}\baseform{2}-\frac{9}{112}\baseform{3}+\frac{1}{56}\baseform{4}+\frac{3}{14}\baseform{5}+\frac{1}{56}\baseform{6}\right)
\\
\bornform{\IrrRep{1750}} =&~ \frac{70}{891}\left(\baseform{1}-\frac{1}{80}\baseform{2}+\frac{1}{112}\baseform{3}-\frac{1}{56}\baseform{4}+\frac{13}{70}\baseform{5}+\frac{1}{70}\baseform{6}\right)
\\
\bornform{\IrrRep{5103}} =&~ \frac{63}{575}\left(\baseform{1}+\frac{3}{112}\baseform{2}+\frac{1}{112}\baseform{3}-\frac{1}{144}\baseform{4}-\frac{1}{42}\baseform{5}-\frac{1}{84}\baseform{6}\right)
\\
\bornform{\IrrRep{5544}} =&~ \frac{56}{225}\left(\baseform{1}-\frac{1}{224}\baseform{2}-\frac{1}{224}\baseform{3}-\frac{1}{224}\baseform{4}-\frac{1}{14}\baseform{5}+\frac{1}{56}\baseform{6}\right)
\\
\bornform{\IrrRep{9625}} =&~ \frac{35}{81}\left(\baseform{1}-\frac{1}{80}\baseform{2}-\frac{1}{560}\baseform{3}+\frac{1}{112}\baseform{4}+\frac{1}{70}\baseform{5}-\frac{1}{140}\baseform{6}\right)
\end{align*}
}
Time \textbf{C}: \Math{1} min \Math{29} sec. Time \textbf{Maple}: \Math{1} sec.
\subsubsection{Representations \Math{\PermRep{22275c}} and \Math{\PermRep{22275d}}}
Rank: \Math{13}.\\ Suborbit lengths: \Math{1, 672, 140, 3360, 5040, 2240, 210, 420, 112, 1680, 1680, 3360, 3360}.
\MathEq{\PermRep{22275c,d}\cong\IrrRep{1}\oplus\IrrRep{22}\oplus\underbrace{\vect{\IrrRep{252}\oplus\IrrRep{252}}}_{}\oplus\underbrace{\vect{\IrrRep{1750}\oplus\IrrRep{1750}}}_{}\oplus\,\IrrRep{3520}\oplus\IrrRep{5103}\oplus\IrrRep{9625}}
\vspace*{-20pt}
{\cmath	 
\begin{align*}
	\bornform{\IrrRep{1}} =&~ \frac{1}{22275}\sum_{k=1}^{13}\baseform{k}
\\
\bornform{\IrrRep{22}} =&~ \frac{2}{2025}\left(\baseform{1}+\frac{13}{28}\baseform{2}-\frac{17}{28}\baseform{3}-\frac{1}{14}\baseform{4}+\frac{11}{56}\baseform{5}+\frac{11}{56}\baseform{6}-\frac{1}{14}\baseform{7}
\right.\\&\left.\hspace*{32pt}
+\frac{13}{28}\baseform{8}+\frac{11}{56}\baseform{9}-\frac{19}{56}\baseform{10}-\frac{19}{56}\baseform{11}-\frac{1}{14}\baseform{12}-\frac{1}{14}\baseform{13}
\right)\\
\bornform{\IrrRep{252}} =&~ \frac{28}{2475}\left\{\baseform{1}-\frac{1}{7616}\vect{129-25\imu\sqrt{671}}\baseform{2}-\frac{3}{476}\vect{1-\imu\frac{5\sqrt{671}}{2}}\baseform{3}
\right.\\&\left.\hspace*{32pt}
-\frac{1}{1088}\vect{\frac{173}{7}+5\imu\sqrt{671}}\baseform{4}
-\frac{5}{7616}\vect{167-\imu\sqrt{671}}\baseform{5}
+\frac{41}{224}\baseform{6}
\right.\\&\left.\hspace*{32pt}
+\frac{1}{952}\vect{257-5\imu\sqrt{671}}\baseform{7}
+\frac{5}{272}\vect{15+\imu\frac{\sqrt{671}}{7}}\baseform{8}
-\frac{3}{7616}\vect{383-25\imu\sqrt{671}}\baseform{9}
\right.\\&\left.\hspace*{32pt}
+\frac{3}{1088}\vect{7+\imu\frac{5\sqrt{671}}{7}}\baseform{10}
+\frac{3}{1088}\vect{7+\imu\frac{5\sqrt{671}}{7}}\baseform{11}
\right\}
\end{align*}\begin{align*}
\bornform{\IrrRep{252}}{'} =&~ \frac{28}{2475}\left\{\baseform{1}
+\frac{1}{7616}\vect{1761-25\imu\sqrt{671}}\baseform{2}
+\frac{3}{476}\vect{\frac{107}{5}-\imu\frac{5\sqrt{671}}{2}}\baseform{3}
\right.\\&\left.\hspace*{32pt}
+\frac{1}{1088}\vect{\frac{3041}{35}+5\imu\sqrt{671}}\baseform{4}+\frac{1}{2828}\vect{\frac{1999}{5}-5\imu\sqrt{671}}\baseform{5}-\frac{173}{1120}\baseform{6}
\right.\\&\left.\hspace*{32pt}
-\frac{1}{952}\vect{121-5\imu\sqrt{671}}\baseform{7}+\frac{1}{272}\vect{\frac{33}{5}-\imu\frac{5\sqrt{671}}{7}}\baseform{8}
\right.\\&\left.\hspace*{32pt}
-\frac{3}{7616}\vect{653+25\imu\sqrt{671}}\baseform{9}+\frac{3}{1088}\vect{\frac{101}{5}-\imu\frac{5\sqrt{671}}{7}}\baseform{10}
\right.\\&\left.\hspace*{32pt}
+\frac{3}{1088}\vect{\frac{101}{5}-\imu\frac{5\sqrt{671}}{7}}\baseform{11}-\frac{3}{40}\baseform{12}-\frac{3}{40}\baseform{13}
\right\}
\\
\bornform{\IrrRep{1750}} =&~ \frac{70}{891}\left\{\baseform{1}-\frac{1}{3472}\vect{185+3\imu\sqrt{55}}\baseform{2}+\frac{3}{1736}\vect{\frac{191}{5}-3\imu\sqrt{55}}\baseform{3}
\right.\\&\left.\hspace*{32pt}
+\frac{1}{3472}\vect{37+\imu\frac{3\sqrt{55}}{5}}\baseform{4}+\frac{3}{2480}\vect{3-\imu\frac{\sqrt{55}}{7}}\baseform{5}-\frac{1}{35}\baseform{6}
\right.\\&\left.\hspace*{32pt}
+\frac{1}{217}\vect{16-\imu\frac{3\sqrt{55}}{2}}\baseform{7}+\frac{1}{1240}\vect{\frac{799}{7}+3\imu\sqrt{55}}\baseform{8}
\right.\\&\left.\hspace*{32pt}
+\frac{3}{3472}\vect{185+3\imu\sqrt{55}}\baseform{9}-\frac{9}{2480}\vect{3-\imu\frac{\sqrt{55}}{7}}\baseform{10}-\frac{9}{2480}\vect{3-\imu\frac{\sqrt{55}}{7}}\baseform{11}
\right\}
\\
\bornform{\IrrRep{1750}}{'} =&~ \frac{70}{891}\left\{\baseform{1}+\frac{1}{3472}\vect{247+3\imu\sqrt{55}}\baseform{2}-\frac{3}{1736}\vect{\frac{377}{5}-3\imu\sqrt{55}}\baseform{3}
\right.\\&\left.\hspace*{32pt}
-\frac{1}{17360}\vect{247+3\imu\sqrt{55}}\baseform{4}-\frac{3}{17360}\vect{83-\imu\sqrt{55}}\baseform{5}-\frac{1}{35}\baseform{6}
\right.\\&\left.\hspace*{32pt}
+\frac{1}{217}\vect{\frac{44}{5}+\imu\frac{3\sqrt{55}}{2}}\baseform{7}+\frac{1}{1240}\vect{1-3\imu\sqrt{55}}\baseform{8}+\frac{3}{3472}\vect{1-3\imu\sqrt{55}}\baseform{9}
\right.\\&\left.\hspace*{32pt}
+\frac{3}{17360}\vect{1-3\imu\sqrt{55}}\baseform{10}
+\frac{3}{17360}\vect{1-3\imu\sqrt{55}}\baseform{11}+\frac{3}{140}\baseform{12}+\frac{3}{140}\baseform{13}
\right\}
\\
\bornform{\IrrRep{3520}} =&~ \frac{64}{405}\left(\baseform{1}
-\frac{1}{224}\baseform{2}+\frac{13}{140}\baseform{3}+\frac{1}{280}\baseform{4}-\frac{1}{280}\baseform{5}-\frac{1}{2240}\baseform{6}-\frac{1}{14}\baseform{7}
\right.\\&\left.\hspace*{32pt}
+\frac{1}{70}\baseform{8}-\frac{13}{112}\baseform{9}-\frac{1}{70}\baseform{10}-\frac{1}{70}\baseform{11}+\frac{11}{1120}\baseform{12}+\frac{11}{1120}\baseform{13}
\right)
\\
\bornform{\IrrRep{5103}} =&~ \frac{63}{275}\left(\baseform{1}+\frac{1}{56}\baseform{2}+\frac{3}{140}\baseform{3}+\frac{1}{105}\baseform{4}-\frac{1}{180}\baseform{5}+\frac{3}{280}\baseform{6}-\frac{1}{42}\baseform{7}
\right.\\&\left.\hspace*{32pt}
-\frac{19}{420}\baseform{8}+\frac{3}{28}\baseform{9}-\frac{1}{420}\baseform{10}-\frac{1}{420}\baseform{11}-\frac{1}{280}\baseform{12}-\frac{1}{280}\baseform{13}
\right)
\\
\bornform{\IrrRep{9625}} =&~ \frac{35}{81}\left(\baseform{1}-\frac{1}{56}\baseform{2}-\frac{1}{28}\baseform{3}-\frac{1}{140}\baseform{4}+\frac{1}{140}\baseform{5}+\frac{1}{280}\baseform{6}+\frac{1}{70}\baseform{7}
\right.\\&\left.\hspace*{32pt}
-\frac{1}{140}\baseform{8}-\frac{1}{28}\baseform{9}+\frac{1}{140}\baseform{10}+\frac{1}{140}\baseform{11}-\frac{1}{280}\baseform{12}-\frac{1}{280}\baseform{13}
\right)
\end{align*}
}
Time \textbf{C}: \Math{3} min \Math{39} sec. Time \textbf{Maple}: \Math{11} sec.
\subsection{Held group \Math{He}}
\paragraph{Main properties:} 
\Math{\cabs{He}=4030387200=2^{10}\cdot3^3\cdot5^2\cdot7^3\cdot17.}~~
\Math{\mathrm{M}\farg{He}=1.}~~ 
\Math{\mathrm{Out}\farg{He}=\CyclG{2}.}
\subsubsection{\Math{29155}-dimensional representation of \Math{He}}
Rank: \Math{12}. Suborbit lengths: \Math{1, 120, 11520, 5760, 1440, 384, 2160, 90, 2880, 2880, 960, 960}.
\MathEq{\PermRep{29155}\cong\IrrRep{1}\oplus\IrrRep{51}\oplus\overline{\IrrRep{51}}\oplus\IrrRep{680}\oplus\underbrace{\vect{\IrrRep{1275}\oplus\IrrRep{1275}}}_{}\oplus\,\IrrRep{1920}\oplus\IrrRep{4352}\oplus\IrrRep{7650}\oplus\IrrRep{11900}}
\vspace*{-20pt}
{\cmath	 
\begin{align*}
	\bornform{\IrrRep{1}} =&~ \frac{1}{29155}\sum_{k=1}^{12}\baseform{k}
\\
\bornform{\IrrRep{51}} =&~ \frac{3}{1715}\left\{\baseform{1}+\frac{5}{12}\baseform{2}-\frac{1}{48}\baseform{3}+\frac{1}{8}\baseform{4}+\frac{1}{8}\baseform{5}\frac{13}{48}\baseform{6}-\frac{1}{6}\baseform{7}-\frac{1}{6}\baseform{8}
\right.\\&\left.\hspace*{32pt}
-\frac{1}{32}\vect{3-\imu\frac{7\sqrt{7}}{7}}\baseform{9}-\frac{1}{32}\vect{3+\imu\frac{7\sqrt{7}}{7}}\baseform{10}
\right.\\&\left.\hspace*{32pt}
+\frac{1}{96}\vect{5+7\imu\sqrt{7}}\baseform{11}+\frac{1}{96}\vect{5-7\imu\sqrt{7}}\baseform{12}
\right\}
\\
\bornform{\IrrRep{680}} =&~ \frac{8}{343}\left(\baseform{1}+\frac{3}{10}\baseform{2}-\frac{1}{48}\baseform{3}-\frac{23}{1440}\baseform{4}-\frac{1}{20}\baseform{5}+\frac{1}{8}\baseform{6}+\frac{1}{120}\baseform{7}
\right.\\&\left.\hspace*{32pt}
+\frac{13}{90}\baseform{8}+\frac{1}{36}\baseform{9}+\frac{1}{36}\baseform{10}+\frac{1}{15}\baseform{11}+\frac{1}{15}\baseform{12}
\right)
\\
\bornform{\IrrRep{1275}} =&~ \frac{15}{343}\left\{\baseform{1}+\frac{1}{4280}\vect{\frac{331}{3}-7\imu\sqrt{231}}\baseform{2}-\frac{1}{25680}\vect{13-\imu\frac{7\sqrt{231}}{3}}\baseform{3}\right.\\&\left.\hspace*{32pt}
-\frac{1}{25680}\vect{\frac{1381}{3}+7\imu\sqrt{231}}\baseform{4}+\frac{1}{25680}\vect{2101+7\imu\sqrt{231}}\baseform{5}
\right.\\&\left.\hspace*{32pt}
-\frac{1}{1712}\vect{13-\imu\frac{7\sqrt{231}}{3}}\baseform{6}+\frac{1}{2568}\vect{\frac{109}{3}-\imu\frac{7\sqrt{231}}{5}}\baseform{7}
\right.\\&\left.\hspace*{32pt}
+\frac{1}{4815}\vect{1571-\imu\frac{7\sqrt{231}}{2}}\baseform{8}-\frac{1}{38520}\vect{467-7\imu\sqrt{231}}\baseform{9}
\right.\\&\left.\hspace*{32pt}
-\frac{1}{38520}\vect{467-7\imu\sqrt{231}}\baseform{10}
\right\}
\\
\bornform{\IrrRep{1275}}{'} =&~ \frac{15}{343}\left\{\baseform{1}+\frac{1}{4280}\vect{\frac{1381}{3}+7\imu\sqrt{231}}\baseform{2}+\frac{1}{25680}\vect{227-\imu\frac{7\sqrt{231}}{3}}\baseform{3}
\right.\\&\left.\hspace*{32pt}
-\frac{1}{25680}\vect{\frac{331}{3}-7\imu\sqrt{231}}\baseform{4}-\frac{1}{25680}\vect{389+7\imu\sqrt{231}}\baseform{5}
\right.\\&\left.\hspace*{32pt}
+\frac{1}{1712}\vect{227-\imu\frac{7\sqrt{231}}{3}}\baseform{6}+\frac{1}{2568}\vect{\frac{319}{3}+\imu\frac{7\sqrt{231}}{5}}\baseform{7}
\right.\\&\left.\hspace*{32pt}
-\frac{1}{4815}\vect{394-\imu\frac{7\sqrt{231}}{2}}\baseform{8}
-\frac{7}{38520}\vect{\frac{157}{2}+\imu\sqrt{231}}\baseform{9}
\right.\\&\left.\hspace*{32pt}
-\frac{7}{38520}\vect{\frac{157}{2}+\imu\sqrt{231}}\baseform{10}
-\frac{1}{16}\baseform{11}-\frac{1}{16}\baseform{12}
\right\}
\end{align*}\begin{align*}
\bornform{\IrrRep{1920}} =&~ \frac{384}{5831}\left(\baseform{1}+\frac{1}{120}\baseform{2}-\frac{7}{384}\baseform{3}+\frac{1}{120}\baseform{4}+\frac{7}{160}\baseform{5}-\frac{7}{384}\baseform{6}+\frac{1}{120}\baseform{7}-\frac{2}{15}\baseform{8}
\right.\\&\left.\hspace*{32pt}
+\frac{5}{192}\baseform{9}+\frac{5}{192}\baseform{10}-\frac{13}{480}\baseform{11}-\frac{13}{480}\baseform{12}
\right)
\\
\bornform{\IrrRep{4352}} =&~ \frac{256}{1715}\left(\baseform{1}+\frac{1}{8}\baseform{2}+\frac{7}{768}\baseform{3}-\frac{5}{576}\baseform{4}-\frac{7}{128}\baseform{6}-\frac{1}{48}\baseform{7}-\frac{1}{18}\baseform{8}+\frac{1}{576}\baseform{9}
\right.\\&\left.\hspace*{32pt}
+\frac{1}{576}\baseform{10}-\frac{1}{192}\baseform{11}-\frac{1}{192}\baseform{12}
\right)
\\
\bornform{\IrrRep{7650}} =&~ \frac{90}{343}\left(\baseform{1}-\frac{1}{20}\baseform{2}+\frac{1}{120}\baseform{4}-\frac{7}{360}\baseform{5}-\frac{1}{90}\baseform{7}+\frac{1}{10}\baseform{8}+\frac{1}{240}\baseform{9}+\frac{1}{240}\baseform{10}
\right.\\&\left.\hspace*{32pt}
-\frac{1}{80}\baseform{11}-\frac{1}{80}\baseform{12}
\right)
\\
\bornform{\IrrRep{11900}} =&~ \frac{20}{49}\left(\baseform{1}-\frac{1}{20}\baseform{2}-\frac{1}{720}\baseform{4}+\frac{1}{120}\baseform{7}
%\right.\\&\left.\hspace*{32pt}
-\frac{1}{18}\baseform{8}-\frac{1}{180}\baseform{9}-\frac{1}{180}\baseform{10}+\frac{1}{60}\baseform{11}+\frac{1}{60}\baseform{12}
\right)
\end{align*}
}
Time \textbf{C}: \Math{5} min \Math{41} sec. Time \textbf{Maple}: \Math{15} sec.
\subsection{Suzuki  group \Math{Suz}}
\paragraph{Main properties:} 
\Math{\cabs{Suz}=448345497600=2^{13}\cdot3^7\cdot5^2\cdot7\cdot11\cdot13.}~~
\Math{\mathrm{M}\farg{Suz}=\CyclG{6}.}\\
\Math{\mathrm{Out}\farg{Suz}=\CyclG{2}.}
\subsubsection{\Math{65520}-dimensional representation of \Math{2.Suz}}
Rank: \Math{10}. Suborbit lengths: \Math{1, 1, 2816, 3960, 2816, 12672, 20736, 20736, 891, 891}.
\MathEq{\PermRep{65520}\cong\IrrRep{1}\oplus\IrrRep{143}\oplus\IrrRep{364_\alpha}\oplus\IrrRep{364_\beta}\oplus\overline{\IrrRep{364_\beta}}\oplus\IrrRep{5940}\oplus\IrrRep{12012}\oplus\IrrRep{14300}\oplus\IrrRep{16016}\oplus\overline{\IrrRep{16016}}}
\vspace*{-20pt}
{\cmath	 
\begin{align*}
	\bornform{\IrrRep{1}} =&~ \frac{1}{65520}\sum_{k=1}^{10}\baseform{k}
\\
\bornform{\IrrRep{143}} =&~ \frac{11}{5040}\left(\baseform{1}+\baseform{2}+\frac{2}{11}\baseform{3}-\frac{1}{11}\baseform{4}+\frac{2}{11}\baseform{5}-\frac{1}{11}\baseform{6}+\frac{3}{11}\baseform{9}+\frac{3}{11}\baseform{10}
\right)
\\
\bornform{\IrrRep{364_\alpha}} =&~ \frac{1}{180}\left(\baseform{1}+\baseform{2}+\frac{1}{16}\baseform{3}+\frac{1}{6}\baseform{4}+\frac{1}{16}\baseform{5}-\frac{1}{24}\baseform{6}-\frac{1}{144}\baseform{7}-\frac{1}{144}\baseform{8}-\frac{1}{9}\baseform{9}-\frac{1}{9}\baseform{10}
\right)
\\
\bornform{\IrrRep{364_\beta}} =&~ \frac{1}{180}\left(\baseform{1}-\baseform{2}-\frac{1}{8}\baseform{3}+\frac{1}{8}\baseform{5}+\imu\frac{\sqrt{3}}{72}\baseform{7}-\imu\frac{\sqrt{3}}{72}\baseform{8}+\imu\frac{\sqrt{3}}{9}\baseform{9}-\imu\frac{\sqrt{3}}{9}\baseform{10}
\right)
\\
\bornform{\IrrRep{5940}} =&~ \frac{33}{364}\left(\baseform{1}+\baseform{2}+\frac{1}{352}\baseform{3}+\frac{1}{66}\baseform{4}+\frac{1}{352}\baseform{5}+\frac{1}{66}\baseform{6}-\frac{7}{864}\baseform{7}-\frac{7}{864}\baseform{8}+\frac{1}{27}\baseform{9}+\frac{1}{27}\baseform{10}
\right)
\\
\bornform{\IrrRep{12012}} =&~ \frac{11}{60}\left(\baseform{1}+\baseform{2}+\frac{1}{88}\baseform{3}-\frac{1}{66}\baseform{4}+\frac{1}{88}\baseform{5}+\frac{1}{264}\baseform{6}-\frac{1}{33}\baseform{9}-\frac{1}{33}\baseform{10}
\right)
\\
\bornform{\IrrRep{14300}} =&~ \frac{55}{252}\left(\baseform{1}+\baseform{2}-\frac{5}{352}\baseform{3}+\frac{1}{330}\baseform{4}-\frac{5}{352}\baseform{5}-\frac{1}{132}\baseform{6}+\frac{1}{288}\baseform{7}+\frac{1}{288}\baseform{8}+\frac{1}{99}\baseform{9}+\frac{1}{99}\baseform{10}
\right)
\\
\bornform{\IrrRep{16016}} =&~ \frac{11}{45}\left(\baseform{1}-\baseform{2}+\frac{1}{352}\baseform{3}-\frac{1}{352}\baseform{5}-\imu\frac{\sqrt{3}}{288}\baseform{7}+\imu\frac{\sqrt{3}}{288}\baseform{8}+\imu\frac{\sqrt{3}}{99}\baseform{9}-\imu\frac{\sqrt{3}}{99}\baseform{10}
\right)
\end{align*}
}
Time \textbf{C}: \Math{28} min \Math{3} sec. Time \textbf{Maple}: \Math{10} sec.

\vspace*{-10pt}
\begin{thebibliography}{9}
%1
\bibitem{Holt}
Holt,~D.~F., Eick,~B., O'Brien,~E.~A.
\textit{Handbook of Computational Group Theory.}\\  Chapman \& Hall/CRC, 2005.
%2
\bibitem{Kornyak1}
Kornyak~V.~V.\\  Quantum models based on finite groups. \textit{J. Phys.: Conf. Ser.} \textbf{965} 012023, 2018.
\href{http://stacks.iop.org/1742-6596/965/i=1/a=012023}{http://stacks.iop.org/1742-6596/965/i=1/a=012023}\\[4pt]
Modeling Quantum Behavior in the Framework of Permutation Groups. \textit{EPJ Web of Conferences} \textbf{173} 01007, 2018.
\href{https://doi.org/10.1051/epjconf/201817301007}{https://doi.org/10.1051/epjconf/201817301007}
%3
\bibitem{Cameron}
Cameron~P.~J.  \emph{Permutation Groups.}~ Cambridge University Press, 1999.
%
\bibitem{atlas}
Wilson, R. A., et al., {\textsc{Atlas} \textit{of finite group representations.}}\\
\href{http://brauer.maths.qmul.ac.uk/Atlas/v3}{http://brauer.maths.qmul.ac.uk/Atlas/v3}.
\end{thebibliography}
\end{document}